\documentclass[
11pt,%
tightenlines,%
twoside,%
onecolumn,%
nofloats,%
nobibnotes,%
nofootinbib,%
superscriptaddress,%
noshowpacs,%
centertags]%
{article}
\usepackage{amsmath,amssymb}

\newtheorem{definition}{Definition}
\newtheorem{theorem}{Theorem}
\newtheorem{fact}{Fact}

\begin{document}

%\titlerunning{Strong degrees of categoricity and weak density} % for running heads
%\authorrunning{Bazhenov at al.} % for running heads
%\authorrunning{First-Author, Second-Author} % for running heads

\title{A notes on the degrees of relative computable categoricity}
% Splitting into lines is performed by the command \\
% The title is written in accordance with the rules of capitalization.

%\author{\firstname{N.~A.}~\surname{Bazhenov}}
%\email[E-mail: ]{bazhenov@math.nsc.ru}
%\affiliation{Sobolev Institute of Mathematics, 4 Acad. Koptyug Ave., Novosibirsk, 630090 Russia}
%\affiliation{Department of Mathematics and Mechanics, Novosibirsk State University, 2 Pirogova St., Novosibirsk, 630090 Russia}

\author{I.~Sh.~Kalimullin}
%\email[E-mail: ]{Iskander.Kalimullin@kpfu.ru} \affiliation{N.I.
%Lobachevskii Institute of Mathematics and Mechanics, Kazan (Volga
%Region) Federal University, 18 Kremlyovskaya str., Kazan, 420008
%Russia}

 % The date of receipt to the editor, i.e. December 06, 2017

\maketitle

\begin{abstract} % You shouldn't use formulas and citations in the abstract.
We  are studying the degrees in which a computable structure is relatively computably categoricity, i.e., computably categorcial among all non-computable copies of the structure. Unlike the  degrees of computable categoricity we can bound the possible degrees of relative computable categoricity by the oracle $\mathbf 0''$. In the case of rigid structures the bound is in fact $\mathbf 0'$. These estimations are precise, in particular we can build a computable structure which is relatively computably categorical only in the degrees above $\mathbf 0''$.
\end{abstract}

%\subclass{03D30, 03D55} % Enter 2010 Mathematics Subject Classification.

%\keywords{computable isomorphism, rigid structure, Turing degrees, degree of categoricity, computable enumerable sets} % Include keywords separeted by comma.

% Text of article starts here.

\section{Introduction}

It is well-known that isomorphic computable structures can have
\emph{different} algorithmic properties. One of the first examples
of this phenomenon was witnessed by Fr{\"o}hlich and
Shepherdson~\cite{FS56}: they constructed two isomorphic computable
fields $F$ and $G$ such that $F$ has a splitting algorithm, but $G$
has no such algorithm. A simpler example concerns copies of the
standard ordering of natural numbers $(\omega;\leq)$. Clearly,
inside the standard presentation of $(\omega;\leq)$, there is an
algorithm which checks whether two given elements $a$ and $b$ are
adjacent. On the other hand, it is not difficult to build a
computable copy $\mathcal{B}$ of $(\omega;\leq)$ such that the
adjacency relation on $\mathcal{B}$ is not computable.

The notion of \emph{computable categoricity} (or
\emph{autostability}), introduced by Mal'tsev~\cite{Mal61,Mal62}, is
intended to capture the structures whose algorithmic behavior is
relatively tame. A computable structure $\mathcal{S}$ is
\emph{computably categorical} if for any computable copy
$\mathcal{A}$ of $\mathcal{S}$, there is a computable isomorphism
$f\colon \mathcal{A} \cong \mathcal{S}$. Informally speaking, if
$\mathcal{A} \cong \mathcal{B}$ are computable copies of a
computably categorical structure, then $\mathcal{A}$ and
$\mathcal{B}$ share the same algorithmic properties.

Computable categoricity has become a cornerstone of the theory,
which studies the algorithmic complexity of isomorphisms. Following
this line of research, Ash~\cite{Ash86a,Ash86b,Ash87} developed the
structural theory of effective categoricity in the levels of the
hyperarithmetical hierarchy. The recent developments of the area
crystallized in the notions of \emph{categoricity spectrum} and
\emph{degree of categoricity}.

\begin{definition}\label{d0}
    Let $\mathbf{d}$ be a Turing degree. A computable structure $\mathcal{A}$ is \emph{$\mathbf{d}$-com\-pu\-tab\-ly categorical} if for every computable copy $\mathcal{B}$ of $\mathcal{A}$, there is a $\mathbf{d}$-com\-pu\-tab\-le isomorphism from $\mathcal{A}$ onto $\mathcal{B}$. The \emph{categoricity spectrum} of $\mathcal{A}$ is the set
    \[
        \textnormal{CatSpec}(\mathcal{A}) = \{ \mathbf{d}\,\colon \mathcal{A} \text{~is~} \mathbf{d} \text{-computably categorical} \}.
    \]
    A Turing degree $\mathbf{d}$ is the \emph{degree of categoricity} of $\mathcal{A}$ if $\mathbf{d}$ is the least degree in the spectrum $\textnormal{CatSpec}(\mathcal{A})$.
\end{definition}

Degrees of categoricity were introduced
by Fokina, Kalimullin, and Miller~\cite{FKM10}.
The article~\cite{FKM10} proves that every 2-c.e. Turing degree $\mathbf{d}$ is the degree of
categoricity for a computable structure  $\mathcal{A}$. For the case when the degree $\mathbb{d}$ is c.e. the construction of $\mathcal{A}$ can be done more easily. In particular, by~\cite{FKM10} every c.e. degree $\mathbf{d}$ is the degree of
categoricity for a computable rigid structure $\mathcal{A}$ (i.e., with trivial automorphism group).

These examples of structures with a degree categoricity has the following nice built-in property: under some fixed oracle the structure $\mathcal{A}$ becomes computably categorical.

\begin{definition}\label{d1}
A structure ${\mathcal A}$ is {\it computably categorical
on a cone } if there is a degree $\mathbf{x}$ such that for every degree ${\bf a}\ge {\bf x}$ and arbitrary ${\bf a}$-computable copy $\mathcal B\cong\mathcal A$ there exists an $\mathbf{a}$-computable isomorphism from $\mathcal B$ onto $\mathcal A$. In  the case when the structure $\mathcal A$ is computable we also say that $\mathcal A$ is {\it relatively $\mathbf x$-computably categorical}.
\end{definition}

\begin{definition}
A Turing degree $\mathbf d$ is $2$-c.e.  if there is a set $D\in\mathbf d$ such that $D=A-B$ for come c.e. sets $A$ and $B$.
\end{definition}
\begin{theorem}\label{th1} (Fokina, Kalimullin, and Miller~\cite{FKM10}).
Every 2-c.e. Turing degree $\mathbf{d}$ is the degree of
categoricity for a computable structure  $\mathcal{A}$ which is relatively $\mathbf d$-computably categorical. This structure can be chosen rigid if $\mathbf{d}$ is c.e.
\end{theorem}

Moreover, all natural examples of computable categorical structures with  are in fact relatively  computably categorical, e.g., the structures $(\mathbb Q,<)$ and $(\mathbb N,+1)$. Analagously, all natural examples of   $\mathbf a$-computable categorical structures, $\mathbf a<\mathbf 0'$, are relatively $\mathbf a$-computably categorcial, while the structure $(\mathbb N,<)$ is $\mathbf 0'$-computably categorical but not computably categorcial on a cone. 

In comparison with computable categoricity the notion of relative categoricity has a nice syntactic characterization:
\begin{theorem}\label{th1.5} (Ash, Knight, Manasse and Slaman \cite{AKMS89}, Chislom~\cite{Chi90}).
A computable structure $\mathcal A$ is relatively $\mathbf d$-computably categorcial iff   there is a tuple $a_1,\dots,a_k\in \mathcal A$ and a $\mathbf d$-c.e. family $\mathcal F$  (so called, Scott family) of  formulae $\varphi(z_1,\dots,z_n)\in\mathcal F$, $n\in\mathbb N$, with parameters $a_1,\dots,a_k$ satisfying the following conditions:
\begin{enumerate}
\item for every tuple $x_1,\dots x_n\in\mathcal A$ there is a  formula $\varphi(z_1,\dots,z_n)\in\mathcal F$ such that
$$
(\mathcal A,a_1,\dots,a_k)\models \varphi(x_1,\dots,x_n);
$$
\item for each formula $\varphi(z_1,\dots,z_n)\in\mathcal F$ and every tuples $x_1,\dots x_n,y_1,\dots,y_n\in\mathcal A$ with
$$(\mathcal A,a_1,\dots,a_k)\models \varphi(x_1,\dots,x_n),
(\mathcal A,a_1,\dots,a_k)\models \varphi(y_1,\dots,y_n)$$
we have an automorphism $p$ of $(\mathcal A,a_1,\dots,a_k)$ (i.e., an automorphism of $\mathcal A$ fixing the parameters: $p(a_i)=a_i$, $1\leq i\leq n$) such that
$
p(x_1)=y_1,\dots, p(x_n)=y_n.
$
\end{enumerate}
\end{theorem}

The left-to-right implication follows from a forcing construction, but the right-to-left implication follows immediately.

Indeed, let $\mathcal B\cong\mathcal A$. Choose a tuple $b_1,\dots b_k\in\mathcal B$ such that $(\mathcal B,b_1,\dots,b_k)\cong (\mathcal A,a_1,\dots,a_k)$. We can define a $\mathbf d$-computable isomorphism  from $\mathcal A$ onto $\mathcal B$ by the standard back-and-forth induction on $n>k$:

\begin{quote}
  If $n>k$ is even let $a_n$ be the least element in $\mathcal A$ such that $a_n\ne a_i, i<n$. Find a formula
  $\varphi(z_1,\dots,z_n)\in\mathcal F$ and an element $b_n\in \mathcal B$,  $b_n\ne b_i, i<n$, such that
$$(\mathcal A,a_1,\dots,a_k)\models \varphi(a_1,\dots,a_n),
(\mathcal A,a_1,\dots,a_k)\models \varphi(b_1,\dots,b_n).$$

If $n>k$ is odd let $b_n$ be the least element in $\mathcal B$ such that $b_n\ne b_i, i<n$. Find a formula
  $\varphi(z_1,\dots,z_n)\in\mathcal F$ and an element $a_n\in \mathcal B$,  $a_n\ne a_i, i<n$ such that
$$(\mathcal A,a_1,\dots,a_k)\models \varphi(a_1,\dots,a_n),
(\mathcal A,a_1,\dots,a_k)\models \varphi(b_1,\dots,b_n).$$

\end{quote}
Then the correspondence $a_n\mapsto b_n$ is a $\mathbf d$-computable isomorphism  from $\mathcal A$ onto $\mathcal B$.

Definitions~\ref{d0} and \ref{d1} lead to the new notion of relative categoricity spectrum and the degree of relative categorcicity.
\begin{definition}
   The \emph{relative categoricity spectrum} of $\mathcal{A}$ is the set
    \[
        \textnormal{RelCatSpec}(\mathcal{A}) = \{ \mathbf{d}\,\colon \mathcal{A} \text{~is~relatively } \mathbf{d} \text{-computably categorical} \}.
    \]
    A Turing degree $\mathbf{d}$ is the \emph{degree of relative computable categoricity} of $\mathcal{A}$ if $\mathbf{d}$ is the least degree in the spectrum $\textnormal{RelCatSpec}(\mathcal{A})$.
\end{definition}
If  a structure $\mathcal A$ has a degree of relative computable categoricity $\mathbf d$ then for every family $\mathcal F$ from Theorem~\ref{th1.5} there is a Turing operator $\Xi$ such that
$$D=\Xi(\mathcal F), \mbox{where } D\in\mathbf d.$$

In fact, by Theorem~\ref{th1} every $2$-c.e. degree is realized as a degree of relative computable categoricity of some structure. Moreover, the Turing operator $\Xi$ from above does not depend on $\mathcal F$. Here we have the uniform case of the degree of relative computable categoricity:
\begin{definition}
     A Turing degree $\mathbf{d}$ is the \emph{uniform degree of relative computable categoricity} of $\mathcal{A}$ if $\mathcal{A}$ is relatively  $\mathbf{d}$-computably categorical and 
there exists a  Turing operator $\Xi$ and a set $D\in \mathbf d$ such that $D=\Xi(\mathcal F)$ for every Scott family $\mathcal F$ for $\mathcal A.$
\end{definition}
%The paper devoted to the study of possible degrees of relative computable categoricity. The special case is the case of rigid structures where we can locate such degrees only below $\mathbf 0'$.

\section{Degrees of relative computable categoricity of rigid structures}

Note that for rigid structures, i.e., the structures with the trivial automorphism group, the conditions of Theorem~\ref{th1.5} sound easier. Moreover, we one can observe the following 

\begin{fact}\label{th1.6}
If a rigid computable structure $\mathcal A$ is  computably categorical on a cone then $\mathcal A$ is relatively $\mathbf 0'$-computably categorical.
\end{fact}

{\it Proof.} By Theorem~\ref{th1.5} if a rigid structure $\mathcal A$ is  computably categorical on a cone then for some
$a_1,\dots,a_k\in \mathcal A$ and a family $\mathcal F$ of  formulae $\varphi(z_1,\dots,z_n)\in\mathcal F$, $n\in\mathbb N$, with parameters $a_1,\dots,a_k$ satisfying the conditions:
\begin{enumerate}
\item for every tuple $x_1,\dots x_n\in\mathcal A$ there is a  formula $\varphi(z_1,\dots,z_n)\in\mathcal F$ such that
$$
(\mathcal A,a_1,\dots,a_k)\models \varphi(x_1,\dots,x_n);
$$
\item for each formula $\varphi(z_1,\dots,z_n)\in\mathcal F$ and every tuples $x_1,\dots x_n,y_1,\dots,y_n\in\mathcal A$ with
$$(\mathcal A,a_1,\dots,a_k)\models \varphi(x_1,\dots,x_n),
(\mathcal A,a_1,\dots,a_k)\models \varphi(y_1,\dots,y_n)$$
we have 
$
x_1=y_1,\dots, x_n=y_n.
$
\end{enumerate}
Then it is easy too see that $x_1,\dots x_n\in\mathcal A$ we can $\mathbf 0'$-computably  find an existential formula $\varphi_{x_1,\dots x_n}(z_1,\dots,z_n)$  such that
$$
(\mathcal A,a_1,\dots,a_k)\models \varphi_{x_1,\dots,x_n}(x_1,\dots,x_n)
$$
 and

\

$(\forall y_1,\dots,y_n\in\mathcal A)[(\mathcal A,a_1,\dots,a_k)\models \varphi_{x_1,\dots,x_n}(y_1,\dots,y_n)\rightarrow $ 

\hspace{15em}$x_1=y_1\;\&\cdots\&\; x_n=y_n].$

\

Then the $\mathbf 0'$-c.e. family $\mathcal G=\{\varphi_{x_1,\dots,x_n}(z_1,\dots,z_n)\mid x_1,\dots x_n\in\mathcal A\}$ again satisfies the conditions 1 and 2, and so by  Theorem~\ref{th1.5} the structure $\mathcal A$ relatively $\mathbf 0'$-computably categorical.

$\Box$

By Theorem~\ref{th1} every c.e. degree is realized as a degree of relative computable categoricity of a rigid structure. In particular, this means that the estimation in Fact~\ref{th1.6} is presice. The next theorem shows that we can not essentially extend Theorem~\ref{th1} among rigid structures.
\begin{theorem}\label{th2} 
If a rigid structure $\mathcal A$ has a uniform degree of
relative computable categoricity ${\bf d}$ then the degree ${\bf d}$ is c.e.
\end{theorem}

{\it Scetch of proof.} By the previous theorem any set $D\in\mathbf d$ is $\Delta^0_2,$ so we can fix such a $D$ with a computable approximation $D(x)=\lim_sd(x,s)$
Let $\Xi$ be  a Turing operator such that $D=\Xi(\mathcal G)$ for every possible Scott family $\mathcal G$ of $\mathcal A$.  Also fix a $D$-c.e. Scott  family of formulae $\mathcal F$  for the structure $\mathcal A$ with  parameters $a_1,\dots, a_k$.

The Turing operator $\Xi(\mathcal G)$ can be represented via a c.e. set $V$  such that 
$$D(x)=y\iff (\exists u,v)[\langle x,y,u,v\rangle \in V\;\&\;D_u\subseteq\mathcal G\;\&\;D_v\subseteq\overline{\mathcal G}],$$
where $\{D_n\}_{n\in\omega}$ be the canonical numbering of all finite subsets of $\omega$. Fix a computable enumeration  $\{V_p\}_{p\in\omega}$ for $V=\cup_pV_p$.

Below we will say that an existential formula
 $\varphi$ with the parameters $a_1,\dots, a_k$ and free variables $z_1,\dots,z_n$ is  {\em rejected by $\mathcal A\restriction s$} if for some  tuples $(x_1,\dots x_n)\ne (y_1,\dots,y_n)$ we have
$$x_1<s,\dots, x_n<s,y_1<s,\dots,y_n<s$$
and 
$$(\mathcal A\restriction s,a_1,\dots,a_k)\models \varphi(x_1,\dots,x_n),
(\mathcal A\restriction ,a_1,\dots,a_k)\models \varphi(y_1,\dots,y_n),$$
i.e., if the condition 2 fails in $(\mathcal A,a_1,\dots,a_k)$ on elements lesser than $s$. 

Define the c.e. set $E=\cup_t E_t$, where

\medskip

$
E_t=\{\langle x,y,s\rangle < t \mid (\forall p<t)(\forall u<p)(\forall v<p)(\exists \varphi\in D_u)[d(x,p)=y\;\&
$

$\langle x,y,u,v\rangle \in V_p\rightarrow\varphi\mbox{ is rejected by }\mathcal A\restriction t]\}.
$

\medskip

Since $D=\Xi(\mathcal F)$ for every $x\in\omega$ there is an $s$ such that $\langle x,D(x),s\rangle\notin E$. Conversely, if  $\langle x,y,s\rangle\notin E$ then $y=D(x)$ since otherwise for some $\langle x,y,u,v\rangle \in V_p$ the formulae from $D_u$ can not be rejected, and so we will have $D\ne\Xi (\widetilde{\mathcal F})$ for a Scott family $\widetilde{\mathcal F}$ which can be obtained after adding the formulae from $D_u$ into $\mathcal F$ and replacement of formulae from $D_v$ by their longer equivalent versions such that  $D_v\cap\widetilde{\mathcal F}=\emptyset$. Thus, we have $D\le_T E.$

It remains to  show $E\le_T D$. From above if $y\ne D(x)$ then $\langle x,y,s\rangle\in E.$ Suppose now we are deciding whether $\langle x,y,s\rangle\in E$ for the case $y=D(x)$. 
Note that  $D$-computably enumerating  $\mathcal F$ and assuming $D=\Xi(\mathcal F)$   we can find  $u,v$ and $p$ such that $u<p ,v<p$, $d(x,p)=y$,
 $D_u\subseteq\mathcal F$ and $\langle x,y,u,v\rangle \in V_p$. Then
$$\langle x,y,s\rangle\in E\iff \langle x,y,s\rangle\in \cup_{t\le p} E_t,$$
since  the non-rejectable formulae from $D_u$ now prevent $\langle x,y,s\rangle\in E_t$ with $t>p$.

$\Box$

\section{Degrees of relative computable categoricity in general case}.

A generalization of Fact~\ref{th1.6} can be stated as following:

\begin{fact}\label{th3}
If a computable structure $\mathcal A$  is computably categorical
on a cone then the structure is relatively $\mathbf{0''}$-computably categorical. Hence, if in addition $\mathcal A$ has a degree of relative computable categoricity ${\bf d}$, then $\mathbf d\leq \mathbf 0''$.
\end{fact}
 
{\it Proof.} (see the proof of Fact~1.4 in \cite{DHM}).
Let $\mathcal A$ be computable and computably categorical
on a cone. By the result of Ash, Knight, Manasse and Slaman \cite{AKMS89} 
 there is a tuple $a_1,\dots,a_k\in \mathcal A$ and a family of formulae $\mathcal F$ satisfying the conditions 1 and 2 of Theorem~\ref{th1.5}:

 \begin{enumerate}
\item for every tuple $x_1,\dots x_n\in\mathcal A$ there is a  formula $\varphi(z_1,\dots,z_n)\in\mathcal F$ such that
$$
(\mathcal A,a_1,\dots,a_k)\models \varphi(x_1,\dots,x_n);
$$
\item for each formula $\varphi(z_1,\dots,z_n)\in\mathcal F$ and every tuples $x_1,\dots x_n,y_1,\dots,y_n\in\mathcal A$ with
$$(\mathcal A,a_1,\dots,a_k)\models \varphi(x_1,\dots,x_n),
(\mathcal A,a_1,\dots,a_k)\models \varphi(y_1,\dots,y_n)$$
we have an automorphism $p$ of $(\mathcal A,a_1,\dots,a_k)$ such that
$$
p(x_1)=y_1,\dots, p(x_n)=y_n.
$$
In this case we will write for simplicity $(x_1,\dots,x_n)\cong_{aut'}(y_1,\dots y_n)$.
\end{enumerate}

 In the general case, for a computable structure this relation $$(x_1,\dots,x_n)\cong_{aut'}(y_1,\dots y_n)$$ is $\Pi^0_2$ due

\

  $
 (x_1,\dots,x_n)\cong_{aut'}(y_1,\dots y_n)\iff (\forall\mbox{ existential }\Psi)
 $

 \hspace{2em}$[(\mathcal A,a_1,\dots,a_k)\models\Psi(x_1,\dots,x_n)\leftrightarrow (\mathcal A,a_1,\dots,a_k)\models \Psi(y_1,\dots,y_n)].$

\

Hence, using the oracle $\mathbf 0''$ we can compute for every    tuple $x_1,\dots,x_n\in \mathcal A$ the existential formula $\varphi_{x_1,\dots,x_n}(z_1,\dots,z_n)$ as the first existential formula $\varphi(z_1,\dots,z_n)$ satisfying the following $\Pi^0_2$ condition:

\

$
(\forall y_1,\dots,y_n\in \mathcal A)[(\mathcal A,a_1,\dots,a_k)\models \varphi_{x_1,\dots,x_n}(y_1,\dots,y_n)\rightarrow 
$

 \hspace{18em}$(x_1,\dots,x_n)\cong_{aut'}(y_1,\dots y_n)]).$

\

Then the $\mathbf 0''$-c.e. family $\mathcal G=\{\varphi_{x_1,\dots,x_n}(z_1,\dots,z_n)\mid x_1,\dots,x_n\in \mathcal A\}$ satisfies the conditions 1 and 2 of Theorem~\ref{th1.5}, so that $\mathcal A$ is relatively $\mathbf{0''}$-computably categorical.

$\Box$

On another hand we  can extend the Theorem~\ref{th1} on a new class of degrees of categoricity. In addition, we show that  the estimation $\mathbf 0''$ in Fact~\ref{th3} is presice.

\begin{definition}
A Turing degree $\mathbf d$ is 2-CEA  if there is a set $D\in\mathbf d$ which is c.e. in a c.e. set $C\le_T D$.
\end{definition}
It is well know that every 2-c.e. degree is in fact 2-CEA: if $D=A-B$ and $A$, $B$ are c.e. then we can built the ``pullback''  $C=f^{-1}(A\cap B)$, where $f$ is a computable function with the range $A$. On another and, the natural example of 2-CEA degree is $\mathbf 0''$ which is clearly not 2-c.e.
\begin{theorem}\label{th4} 
If $\mathbf d\geq \mathbf c$, $\mathbf d$ is c.e. in $ \mathbf c$, and $\mathbf c$ is c.e. then there exists a computable structure  $\mathcal{A}$   which is computably categorical on a cone, and has the degree of categoricit $\mathbf d$.
\end{theorem}

{\it Proof.}
Let  $D\in\mathbf d$ which is c.e. in a c.e. set $C\le_T D$. Then the set $D$ can be be presented in the $\Sigma_2^0$-form
$$x\in D\iff (\exists y)(\forall z\geq y)R(x,z),$$
where the predicate $R$ is computable and the $\Pi^0_1$-relation
$$Q(x,y)\iff (\forall z\geq y)R(x,z) $$
is computable in $C\le_T D$. Introduce the modulus function
$$
m(x)=\begin{cases}\min\{y:Q(x,y)\},&\mbox{if }x\in D;\cr
\infty,& \mbox{if }x\notin D.
\end{cases}
$$
which is computable in the oracle $D$.

Let $\mathcal A$ be the graph with the vertices 
$$u_x \mbox{ for } x\in\mathbb N; v_{x,y} \mbox{ for } x\in\mathbb N, y\leq m(x);   w_{x,y} \mbox{ for } x\in\mathbb N, y<m(x); $$
$$u_x', v_{x,0}' \mbox{ for } x\in\mathbb N; v_{x,1}', w_{x,0}' \mbox{ for } x\in C,$$
and the edge relations of the one of the forms:
$$
\{u_x,u_x'\}, \{u_x',u_{x+1}\},
$$
$$
\{u_x,v_{x,0}\},  \{v_{x,y},v_{x,y+1}\}, \{u_x,w_{x,0}\}, \{w_{x,y},w_{x,y+1}\},
$$
$$
\{u_x',v_{x,0}'\},  \{v_{x,y}',v_{x,y+1}'\}, \{u_x',w_{x,0}'\}, \{w_{x,y}',w_{x,y+1}'\}.
$$
Since the relations ``$x\in C$''  and ``$y<m(x)$'' are c.e. we can effectively map $\mathcal A$ onto a computable graph with the set of vertices equal to $\mathbb N.$

To show that the graph $\mathcal A$ is relatively $\mathbf d$-computably categorical suppose that $\mathcal B\cong \mathcal A$ is $\mathbf x$-computable for $\mathbf x\ge \mathbf d.$ If we fix in $\mathcal B$ the element corresponding to $u_0$ we then can effectively find inn $\mathcal B$ other elements of the chain $u_0-u_0'-u_1-u_1'-u_2-\cdots$. Since $\mathbf x$ knows the membership in $C$ we can $\mathbf x$-computably map the elements $v_{x,0}'$ and, if any,   $v_{x,1}'$ and $w_{x_0}'$. Since $m(x)$ is  $\mathbf x$-computable too we can effectively map also the elements of the chains $v_{x,0}-v_{x,1}-\cdots$ and $w_{x,0}-w_{x,1}-\cdots$. Thus, $\mathcal A$ is relatively $\mathbf d$-computably categorical . Moreover, one can find precisely a $\mathbf d$-c.e. family of formulae $\mathcal F$ in the terms of Theorem~\ref{th1.5} with the parameter $a_1=u_0.$

To show that the degree $\mathbf d$ is the degree of relative computable categoricity suppose that  $\mathcal A$ is relatively $\mathbf x$-computably categorical. Then for some tuple $a_1,\dots,a_k\in\mathcal A$ a $\mathbf d$-c.e. family of formulae $\mathcal F$ satisfying the conditions 1 and 2 of Theorem~\ref{th1.5}. Let $x_0$ be large enough such that for $x\ge x_0$ the chains $v_{x,0}-v_{x,1}-\cdots$, $w_{x,0}-w_{x,1}-\cdots$, $v_{x,0}'-v_{x,1}'-\cdots$, $w_{x,0}'-w_{x,1}'-\cdots$ do not contain the elements $a_1,\dots,a_k$.

Let $C=\cup_s C_c$ be the computable enumeration of $C$ via finite approxiimations $C_s, s\in\mathbb N$,
$$Q(x,y,s)\iff (\forall z)[y\ge z<s\rightarrow R(x,z)],$$
$$
m(x,s)=\begin{cases}\min\{y:Q(x,y,s)\},&\mbox{if }Q(x,y,s)\mbox{for some }y<s;\cr
s& \mbox{otherwise,}
\end{cases}
$$
and finally let $\mathcal A_s$ be the finite subgraph of $\mathcal A$ on the vertices
$$u_x \mbox{ for } x<s; v_{x,y} \mbox{ for } x<s, y\leq m(x,s);   w_{x,y} \mbox{ for } x<s, y<m(x,s); $$
$$u_x', v_{x,0}' \mbox{ for } x<s; v_{x,1}', w_{x,0}' \mbox{ for } x\in C_s.$$
Now to $\mathbf x$-computably decide the membership  $x\in C$, $x\ge x_0$, it is enough to find a formula $\varphi(x)\in\mathcal F$ such that
$(\mathcal A_s,a_1,\dots,a_k)\models \varphi(v_{x,0}')$
for some $s.$ Then 
$$
x\in C\iff x\in C_s,
$$  
since otherwise the formula $\varphi$ is not able to distinguish the elements $v_{x,0}'$ and $w_{x,0}'$. Thus, the set $C$ is $\mathbf x$-computable, and hence $D$ is $\mathbf x$-c.e.

Note that the elements  $v_{x,0}$ and $w_{x,0}$ are in the same automorphism orbit iff $x\notin D$. Therefore, for every $x\ge x_0$
$$
x\notin D\iff (\exists s)(\exists \varphi(x)\in\mathcal F)[(\mathcal A_s,a_1,\dots,a_k)\models \varphi(v_{x,0})\;\&\; \varphi(w_{x,0})], 
$$
so that $D$ is also co-$\mathbf x$-c.e., and hence $\mathbf d\le\mathbf x.$

$\Box$
%

%there is a nontrivial initial interval $\mathcal{I}$ in Turing
%degrees such that any non-c.e. degree from $\mathcal{I}$ cannot be a
%degree of categoricity for a rigid structure. More precisely, the
%following result is obtained:

%This section is devoted to the proof of result which is a base for
%Theorem \ref{maintheor}. In \cite{BKY16} it was shown that a strong
%degree of categority is equivalent to being a degree of a structure
%with spectral dimension one. Thus, Theorem \ref{helptheor} together
%with Lemma \ref{tospecdim1} lead to Theorem \ref{maintheor}.

%\begin{theorem} \label{th2}
%There exists a noncomputable c.e. set $C$ such that for
%any $D \leq_T C$, one of the following holds: either $\deg(D)$ is
%c.e., or $\deg(D)$ is not a degree of categoricity.
%\end{theorem}

%%%%%%%%%%%%%%%%%%%%%%%%%%%%

%\begin{acknowledgments}
   % The work was supported by  %the Russian Science Foundation, project No.
%18-11-00028. %Bazhenov was supported by the grant ....  Kalimullin
%was supported by ... . Yamaleev was supported by the Russian Science
%Foundation, project No. 18-11-00028.

%\end{acknowledgments}

%%%%%%%%%%%%%%%%%%%%%

\end{document}